\documentclass[12pt,a4paper]{amsart}
\usepackage{amssymb}
\usepackage{mathrsfs}

\headsep=1cm \numberwithin{equation}{section}
\newtheorem{thm}{Theorem}[section]
\newtheorem{cor}[thm]{Corollary}
\newtheorem{lem}[thm]{Lemma}
\newtheorem{prop}[thm]{Proposition}
\theoremstyle{definition}
\newtheorem{defn}[thm]{Definition}
\theoremstyle{remark}

\numberwithin{equation}{section}

\newcommand{\To}{\longrightarrow}

\newcommand{\nat}{\mathbb N}
\newcommand{\m}{\frak m }
\newcommand{\p}{\frak p }

\newcommand\Supp{\operatorname{Supp}}
\newcommand\Ass{\operatorname{Ass}}

\newcommand\Spec{\operatorname{Spec}}

\newcommand\Hom{\operatorname{Hom}}
\newcommand\Ext{\operatorname{Ext}}
\newcommand\grade{\operatorname{grade}}

\newcommand\Max{\operatorname{Max}}
\newcommand\depth{\operatorname{depth}}
\begin{document}
\title[Faltings' local-global principle of local cohomology modules]
{Faltings' local-global principle for the in dimension $\bf< n$  of local cohomology modules}
\author[R.  Naghipour, R.  Maddahali and Kh.  Ahmadi Amoli]{ Reza Naghipour$^*$,  Robabeh Maddahali and Khadijeh Ahmadi Amoli\\\\\\\,
\vspace*{0.5cm} Dedicated to Professor Peter Schenzel}
\address{Department of Mathematics, University of Tabriz, Tabriz, Iran;
and School of Mathematics, Institute for Research in Fundamental
Sciences (IPM), P.O. Box 19395-5746, Tehran, Iran.}
\email{naghipour@ipm.ir} \email {naghipour@tabrizu.ac.ir}
\address{Department of Mathematics, University of Payame Noor, P. O.  Box 19395-3697, Tehran, Iran.}
\email{maddahali@phd.pnu.ac.ir}
\address{Department of Mathematics, University of Payame Noor, P. O.  Box 19395-3697, Tehran, Iran.}
\email{khahmadi@pnu.ac.ir}

\thanks{ 
This research  was  in part supported by a grant from IPM.\\
$^*$Corresponding author: e-mail: {\it naghipour@ipm.ir} (Reza Naghipour)}%

\begin{abstract}
The concept of Faltings' local-global principle for the  in dimension $< n$ of local cohomology modules
over a  Noetherian ring $R$ is introduced, and it is shown that this principle holds at levels 1, 2. We also
establish the same principle at all levels over an arbitrary  Noetherian ring of dimension not exceeding 3.
These generalize  the main results of Brodmann et al. in \cite{BRS}. Moreover, as a generalization of Raghavan's result, we show that the Faltings' local-global principle for the
in dimension $<n$ of local cohomology modules  holds at all levels $r\in \mathbb{N}$  whenever the ring $R$ is a homomorphic image of a Noetherian Gorenstein ring.
Finally, it is shown that if $M$ is a
finitely generated $R$-module, $\frak a$ an ideal of $R$  and $r$ a non-negative integer such that $\frak a^tH^i_{\frak a}(M)$
is  in dimension $< 2$ for all $i<r$ and for some positive integer $t$, then for any minimax submodule $N$ of  $H^r_{\frak a}(M)$,
the $R$-module $\Hom_R(R/\frak a,  H^r_{\frak a}(M)/N)$ is finitely generated. As a consequence, it follows that the associated
primes of $H^r_{\frak a}(M)/N$ are finite.  This generalizes the main results of Brodmann-Lashgari \cite{BL} and Quy \cite{Qu}.
\end{abstract}
\maketitle
\section{Introduction}
Let  $R$ denote a commutative Noetherian ring
(with non-zero identity) and $\frak a$ an ideal of $R$. For an $R$-module $M$, the
$i{\rm th}$ local cohomology module of $M$ with support in $V(\frak a)$
is defined as:
$$H^i_{\frak a}(M) = \underset{n\geq1} {\varinjlim}\,\, {\rm Ext}^i_R(R/\frak a^n, M).$$
Local cohomology was defined and studied by Grothendieck. We refer the reader to \cite{BS} or \cite{Gr1} for more
details about local cohomology.  An important theorem in local
cohomology is Faltings' local-global principle for the finiteness
dimension of local cohomology modules \cite[Satz 1]{Fa1}, which
states that for a positive integer $r$,   the $R_{\frak p}$-module $H^i_{\frak aR_{\frak p}}(M_{\frak p})$ is finitely generated
for all $i\leq r$ and for all ${\frak p}\in  \Spec(R)$ if and only if
the $R$-module $H^i_{\frak a}(M)$ is finitely generated for all $i\leq r$.\\
Another formulation of Faltings' local-global principle,
particularly relevant for this paper, is in terms of the
generalization of the finiteness dimension $f_{\frak a}(M)$ of $M$ relative
to $\frak a$, where
$$f_{\frak a}(M):=\inf\{i\in \Bbb{N}_0\,\,|\,\,H^i_{\frak a}(M)\,\,{\rm is}\,\,{\rm not}\,\,{\rm finitely}\,\,{\rm generated}
\};\,\,\,\,\,\,\,\,\,\,\,\,\,\,\,\,\,\,\,\,\,\,\,\,\,(\dag)$$
 with the usual convention that the
infimum of the empty set of integers is interpreted as $\infty$.
For any non-negative integer $n$, the {\it $n{\rm th}$ finiteness dimension} $f^n_{\frak a}(M)$ of $M$
relative to $\frak a$ is defined  by $$f^n_{\frak a}(M):=\inf\{f_{\frak aR_{\frak p}}(M_{\frak
p})\,\,|\,\,{\frak p}\in {\rm Supp}(M/\frak a M)\,\,{\rm and}\,\,  \dim R/{\frak p}\geq n\}.$$ Note that $f^n_{\frak a}(M)$ is either a positive
integer or $\infty$ and that $f^0_{\frak a}(M)=f_{\frak a}(M)$. The {\it $n{\rm th}$ finiteness dimension} $f^n_{\frak a}(M)$ of $M$
relative to $\frak a$ has been introduced by Bahmanpour et al. in \cite{BNS1}.

Recall that the $\frak b$-{\it finiteness dimension} of $M$ relative to $\frak a$ is defined by
\begin{eqnarray*}
f_{\frak a}^{\frak b}(M)&:=& \inf\{i\in \Bbb{N}_0\,\,|\,\,\frak b\not\subseteq {\rm Rad}(0:_{R}H^i_{\frak a}(M))\}\\
&=& \inf\{i\in \Bbb{N}_0\,\,|\,\,{\frak b}^nH^i_{\frak a}(M)\neq 0\,\,{\rm for}\,\,{\rm all}\,\,n\in \Bbb{N}\},
\end{eqnarray*}
where $\frak b$ is a second ideal of $R$.

 Brodmann et al. in \cite{BRS} defined and studied the concept of the local-global principle for annihilation
 of local cohomology modules at level $r\in\mathbb{N}$ for the ideals $\frak a$ and $\frak b$ of $R$. We say that the local-global principle for the annihilation of local cohomology modules
 holds at level $r$ if for every choice of ideals ${\frak a}$, ${\frak b}$ of $R$ and every choice of finitely generated $R$-module $M$,
 it is the case that  $$f_{{\frak a}R_{\frak p}}^{{\frak b}R_{\frak p}}(M_{\frak p})>r \,\,\,\,\, \text{ for all } {\frak p}\in\Spec(R) \Longleftrightarrow f_{\frak a}^{\frak b}(M)>r.$$
It is shown in \cite{BRS} that the local-global principle for the annihilation of local cohomology modules holds at levels 1, 2, over an arbitrary commutative Noetherian ring $R$ and at all
 levels whenever $\dim R\leq 4$.

For a non-negative integer $n$ we say that $M$ is {\it in dimension} $<n$, if $\dim\Supp (M/N)< n$ for some finitely generated submodule $ N $ of $ M$.  In \cite{AN1}, Asdollahi and Naghipour
  introduced the notion   $h_{\frak a}^n(M)$ as follows:
$$h_{\frak a}^n(M)={\rm inf} \{ i\in\mathbb{N}_0: \,\, H^{i}_{\frak a}(M) \,\, \text{is not in dimension}< n\}. $$
This motivates to introduce the notion of  $h_{\frak a}^{\frak b}(M)^n$   by
           $${h_{\frak a}^{\frak b}(M)}^{n}:=\inf\{i\in\mathbb{N}_0:{\frak
        b}^tH_{\frak a}^i(M)\text{ is not in dimension $<n$  for all}\,\, t \in\mathbb{N}\}.$$

         Note that,  $h_{\frak a}^{\frak b}(M)^n$ is  either a non-negative integer or $\infty$, and if $M$ is a finitely generated $R$-module then
          $h_{\frak a}^{\frak a} (M)^n=h_{\frak a}^n(M)$  and that ${h_{\frak a}^{\frak b} (M)}^0={ f_{\frak a}^{\frak b} (M)}$.

We say that the local-global principle for the in dimension $<n$  of local cohomology modules holds at level
$r\in\mathbb{N}$ if for every choice of ideals ${\frak a}$, ${\frak b}$ of
$R$ with ${\frak b}\subseteq {\frak a}$ and every choice of finitely generated $R$-module $M$, it is the case that
$$h_{{\frak a}R_{\frak p}}^{{\frak b}R_{\frak p}}(M_{\frak p})^n>r\,\,\,\,\, \text{ for all }
{\frak p}\in\Spec(R) \Longleftrightarrow h_{{\frak a}}^{{\frak b}}(M)^n>r.$$

Our main result in Section 2 is to introduce the concept of Faltings' local-global principle for the in dimension $<n$ of local cohomology modules
over a commutative Noetherian ring $R$, and we show that this principle holds at levels 1, 2. We also
establish the same principle at all levels over an arbitrary commutative Noetherian ring of dimension not exceeding 3.
Our tools for proving the main result in Section 2 is the following:
\begin{thm}
Suppose that $R$ is a  Noetherian ring  and let  $\frak a, {\frak b}$ be  two ideals of $R$ such that ${\frak b}\subseteq{\frak a}$.  Assume that  $M$ is a
finitely generated $R$-module and let  $r$ be a positive integer such that  the local cohomology
modules $H_{\frak a}^0(M),\dots, H_{\frak a}^{r-1}(M)$ are ${\frak a}$-cofinite. Then
$$h_{{\frak a}R_{\p}}^{{\frak b}R_{\p}}(M_{\p})^n>r \text{ for all } {\p}\in \Spec (R) \Longleftrightarrow h_{{\frak a}}^{{\frak b}}(M)^n>r.$$
\end{thm}

Pursuing this point of view further we establish the following consequence of Theorem 1.1  which is an extension
of the results of Brodmann et al. in \cite[Corollary 2.3]{BRS} and Raghavan in \cite{Ra} for an arbitrary Noetherian ring.
\begin{cor}
Let $R$ be a  Noetherian ring, $M$ a finitely generated $R$-module and, $\frak a, \frak b$ two ideals of $R$ such that ${\frak b}\subseteq{\frak a}$
and $\frak aM\neq M$. Set  $r\in\{1, \grade_M\frak a, f_{\frak a}(M), f_{\frak a}^1(M), f_{\frak a}^2(M)\}$. Then
$$h_{{\frak a}R_{\p}}^{{\frak b}R_{\p}}(M_{\p})^n>r \,\,\,\,  \text{ for all } {\p}\in \Spec (R) \Longleftrightarrow h_{{\frak a}}^{{\frak b}}(M)^n>r.$$
\end{cor}

Moreover in this section, we explore an interrelation between this principle and the Faltings' local-global principle for the annihilation of local
cohomology modules, and show that the local-global principle for the annihilation of local cohomology modules holds at levels  $1, 2$ over $R$ and at all levels whenever $\dim R\leq 3$.
These generalize and reprove the main results of Brodmann et al. in \cite{BRS}.

In \cite{Ra}, Raghavan deduced from the Faltings' Theorem for the annihilation of local cohomology modules \cite{F3} that  if $R$ is a homomorphic image of a Noetherian regular ring,
then the local-global principle holds at all levels $r\in \mathbb{N}$.  In  Section 3, as a generalization of the Raghavan's result, we show that the Faltings' local-global principle for the
in dimension $<n$ of local cohomology modules  holds at all levels $r\in \mathbb{N}$,  whenever the ring $R$ is a homomorphic image of a Noetherian Gorenstein ring.  More precisely we shall
 prove the following:

\begin{thm}
Suppose that $R$ is a  Noetherian ring which  is a homomorphic image of a Gorenstein ring. Then the local-global principle (for the  in dimension $<n$  of local cohomology modules) holds at
all levels $r\in\mathbb{N}$.
\end{thm}

The result in Theorem 1.3 is proved in Theorem 3.1.  Our method is
based on the notion of the $n$th $\frak b$-{\it minimum  $\frak
a$-adjusted depth of $M$} (see \cite{DN2})

\begin{align*}
\lambda_{\frak a}^{\frak b}(M)_n:=\inf\{\lambda_{\frak a R_{\frak p}}^{\frak b
R_{\frak p}}(M_{\frak p})|\, \dim R/ \frak p\geq n\},
\end{align*}
where $\lambda_{\frak a}^{\frak b}(M):=\inf\{ {\depth}\,M_{\frak p}+{{\rm ht}}(\frak a+\frak p/{\frak p)}\mid \frak p\in\Spec (R)\setminus V(\frak b)\}$ is the
  $\frak b$-{\it minimum  $\frak a$-adjusted depth of $M$}, (see \cite[Definition 9.2.2]{BS}).

Also, in this section we show that if   $M$ is a finitely generated $R$-module and  $\frak a,  \frak b$  are ideals of $R$  with ${\frak b}\subseteq{\frak a}$,
then for any  non-negative integer $n$,
$$h_{{\frak a}}^{{\frak b}}(M)^{n}\geq\depth({\frak b},M)   \Longleftrightarrow h_{{\frak a}}^{n}(M)\geq\depth({\frak b},M).$$

Finally, in Section 3 we prove some finiteness results about the associated primes of local cohomology modules. In fact, we will generalize the main results
of Brodmann-Lashgari \cite{BL} and Quy \cite{Qu}.

\begin{thm}
    Let $R$ be a Noetherian ring, $M$  a  finitely generated $R$-module
    and $\frak a$ an ideal of $R$. Let $r$ be a non-negative integer such that $\frak a^tH^{i}_{I}(M)$ is
    in dimension $<2$ for all $i<r$ and for some $t\in\mathbb{N}_0$.  Then, for any minimax submodule $N$ of
    $H^{r}_{\frak a}(M)$,   the $R$-module $\Hom_{R}(R/\frak a, H^{r}_{\frak a}(M)/N)$ is finitely generated. In particular,
    the set $\Ass_{R}(H^{r}_{\frak a}(M)/N)$ is finite.
\end{thm}

Throughout this paper, $R$ will always be a commutative Noetherian
ring with non-zero identity and $\frak a$ will be an ideal of $R$. Recall
that an $R$-module $L$ is called $\frak a$-{\it cofinite} if $\Supp (L)\subseteq V(\frak a)$ and ${\rm Ext}^j_R(R/\frak a, L)$ is finitely
generated for all $j\geq 0$. The concept of $\frak a$-cofinite modules
were introduced by Hartshorne \cite{Ha}.  An $R$-module $L$ is said to be {\it minimax}, if there exists a
finitely generated submodule $N$ of $L$, such that $L/N$ is
Artinian. The class of minimax modules was introduced by H.
Z\"{o}schinger \cite{Zo1} and he has given in \cite{Zo1, Zo2} many
equivalent conditions for a module to be minimax.
We shall use $\Max (R)$ to denote the set of all maximal
ideals of $R$. Also, for any ideal $\frak a$ of $R$, we denote
$\{\frak p \in \Spec(R):\, \frak p\supseteq \frak a \}$ by
$V(\frak a)$. Finally, for any ideal $\frak{b}$ of $R$, the {\it
radical} of $\frak{b}$, denoted by ${\rm Rad}(\frak{b})$, is defined to
be the set $\{x\in R \,: \, x^n \in \frak{b}$ for some $n \in
\mathbb{N}\}$. For any unexplained notation and terminology we refer
the reader to \cite{BS} and \cite{Mat}.

\section{\textbf{Local-global principle  and annihilation of local cohomology modules}}
In this section we investigate the concept of Faltings' local-global principle for the in dimension $<n$ of local cohomology modules
over a commutative Noetherian ring $R$. To this end we begin by the definitions of the {\it in dimension $<n$ modules} and the ${\frak b}$-{\it dimension $<n$ of a module}
relative to an ideal $\frak a$ of $R$.  Then we show the local-global principle for the annihilation of local cohomology modules holds at level $1, 2$ over $R$ and at all levels
 whenever $\dim R\leq 3$. These  extend  the main results of Brodmann et al. in \cite{BRS}.
\begin{defn}\label{def1}
Let $R$ be a  Noetherian ring and $M$ an  $R$-module. For a
non-negative integer $n$ we say that $M$ is in dimension $<n$, if
$\dim\Supp (M/N)< n$ for some finitely generated submodule $N$ of $
M$.
\end{defn}
\begin{defn}\label{def2}
Let $R$ be a  Noetherian ring and $M$ an $R$-module. Let ${\frak a},{\frak b}$ be two ideals of $R$. For a non-negative integer $n$, we define the in ${\frak b}$-dimension $<n$ of
 $M$ relative to ${\frak a}$, denoted by  $h_{\frak a}^{\frak b}(M)^n$,  by
           $${h_{\frak a}^{\frak b}(M)}^{n}:=\inf\{i\in\mathbb{N}_0:{\frak
        b}^tH_{\frak a}^i(M)\text{ is not in dimension $<n$  for all}\,\, t \in\mathbb{N}\}.$$
 Note that,  $h_{\frak a}^{\frak b}(M)^n$ is  either a non-negative integer or $\infty$, and if $M$ is a finitely generated $R$-module in view of \cite[Definition 2.4]{MNS}, \cite [Theorem 2.2]{Ho}
 and Proposition \ref{pro1} that   $h_{\frak a}^{\frak a} (M)^n=h_{\frak a}^n(M)$, and that ${h_{\frak a}^{\frak b} (M)}^0={ f_{\frak a}^{\frak b} (M)}$.
\end{defn}

The following lemma is needed in the proof of Proposition 2.4.

\begin{lem}\label{lem1}
Let $R$ be a  Noetherian ring, ${\frak a}$ an ideal of $R$, and $M$ an arbitrary $R$-module. Then ${\frak
a}M$ is in dimension $<n$ if and only if $M/(0:_M{\frak a})$ is  in dimension $<n$, where $n$ is a non-negative integer.
\end{lem}

\proof  Let ${\frak a}=(a_1, \dots, a_t)$.   Suppose first that ${\frak a}M$ is in dimension $<n$. Then $f:M\longrightarrow (\frak a M)^t$ defined by
$f(x)=(a_1x, \dots, a_tx)$ is an $R$-epimorphism and $Ker(f)=(0:_M{\frak a})$. Therefore,  in view of \cite[Corollary 2.14]{MNS}, $M/(0:_M{\frak a})$ is  in dimension $<n$.
Conversely, if $M/(0:_M{\frak a})$ is  in dimension $<n$ and  $g:M^t\longrightarrow \frak a M$  be the $R$-epimorphism for which
\begin{center}
$g(x_1, \dots, x_t)=\sum_{i=1}^t a_ix_i$,   \,\,\,\,\,\,\,\,\,\,\,\, for all $(x_1, \dots, x_t)\in M^t$,
\end{center}
then $g$ induces  an $R$-epimorphism $g^*:(M/(0:_M{\frak a}))^t\longrightarrow \frak a M$ for which
\begin{center}
$g^*((x_1, \dots, x_t)+(0:_M{\frak a})^t)=g(x_1, \dots, x_t)$, \,\,\,\,\,\,\,\,\,\,\,\, for all $(x_1, \dots, x_t)\in M^t$.
\end{center}
 Now, it follows from \cite[Proposition 2.12 and Corollary 2.14]{MNS} that ${\frak a}M$ is in dimension $<n$.  \qed \\

The following proposition, which is a generalization of \cite[Lemma 9.1.2]{BS}, states that the $R$-modules $H_{\frak a}^0(M)$,\dots, $H_{\frak a}^{s-1}(M)$
are in dimension $<n$ if and only if there is an integer $t\in\mathbb{N}$ such that the $R$-modules ${\frak a}^tH_{\frak a}^0(M)$,\dots, ${\frak a}^tH_{\frak a}^{s-1}(M)$ are in dimension $<n$.
\begin{prop}\label{pro1}
Let $R$ be a  Noetherian ring and ${\frak a}$ an ideal of $R$. Let $s$  and $n$ be two non-negative integers.  Let $M$ be an arbitrary $R$-module such that
 $\Ext^{s-1}_R(R/\frak a, M)$  is in dimension $<n$. Then the following statements are equivalent:

{\rm(i)}  $H_{\frak a}^i(M)$ is in dimension $<n$ for all $i<s$;

{\rm (ii)} There exists  an integer $t\geq1$ such that ${\frak a}^tH_{\frak a}^i(M)$ is in dimension $<n$ for all $i<s$.
\end{prop}

\proof The implication ${\rm(i)}\Longrightarrow {\rm(ii)}$ is obviously true by \cite[Proposition 2.12]{MNS}. In order to show ${\rm(ii)}\Longrightarrow {\rm(i)}$,
we proceed by induction on $s$. If $s=1$, then for some integere $t\geq1$, ${\frak a}^tH_{\frak a}^0(M)$ is in dimension $<n$. Moreover, in view of the assumption, ${\rm Hom_{R}}(R/\frak a, M)$  is in dimension $<n$. Now, since $${\rm Hom_{R}}(R/\frak a,H_{\frak a}^0(M)) \cong  {\rm Hom_{R}}(R/\frak a, M),$$  it follows that the $R$-module $(0:_{H_{\frak a}^{0}(M)}{\frak a})$  is in dimension $<n$. Therefore, it yields from Lemma \ref{lem1} and \cite[Propositions 2.12 and 2.17]{MNS} that  $H_{\frak a}^0(M)$ is in dimension $<n$. Suppose that $s>1$, and the case $s-1$ is settled. By inductive
hypothesis the $R$-module $H_{\frak a}^i(M)$ is in dimension $<n$ for all
$i<s-1$, and so it is enough  to show that the $R$-module $H_{\frak
a}^{s-1}(M)$ is in dimension $<n$. For this purpose, as there is  an integer $t\geq1$  such that ${\frak a}^tH_{\frak a}^{s-1}(M)$ is in dimension $<n$, it follows from
 Lemma \ref{lem1} that $R$-module $H_{\frak
a}^{s-1}(M)/(0:_{H_{\frak a}^{s-1}(M)}{\frak a}^t)$ is in dimension $<n$. On the other hand, by virtue of \cite[Corollary 2.16]{MNS}, the $R$-module
${\rm Ext}^j_R(R/\frak a,H_{\frak a}^i(M))$  is in dimension $<n$ for all $i<s-1$ and all $j\geq0$. Hence, it follows from \cite[Theorem 2.2]{AT} that
 ${\rm Hom}(R/\frak a,H_{\frak a}^{s-1}(M))$ is in dimension $<n$, and so in view of \cite[Proposition 2.17]{MNS},  ${\rm Hom}(R/\frak a^{t},H_{\frak a}^{s-1}(M))$ is also in dimension $<n$.
Consequently, it follows from Lemma \ref{lem1}  and \cite[Corollary 2.13]{MNS} that the $R$-module $H_{\frak a}^{s-1}(M)$ is in dimension $<n$, as required.\qed \\


Before we shall state the next result, we have to recall the notion of the  ${\frak b}$-minimaxness dimension $\mu_{\frak a}^{\frak b}(M)$ of $M$ relative to ${\frak a}$ which is defined by
$$\mu_{\frak a}^{\frak b}(M)=\inf\{i\in\mathbb{N}:{\frak b}^tH_{\frak a}^i(M)\, \text {is not minimax for all
}t\in\mathbb{N}\}.$$
Here $M$ denotes a finitely generated module over a Noetherian ring $R$, and ${\frak a}$, ${\frak b}$ denote two ideals of $R$ such that ${\frak b}\subseteq{\frak a}$
(see \cite[Definition 2.4]{DN1}). Also, for any non-negative integer $n$, Asadollahi and Naghipour in \cite{AN2} defined  the upper  $n$th ${\frak b}$-finiteness dimension $f_{\frak
    a}^{\frak b}(M)^{n}$ of $M$ relative to ${\frak a}$ by
 $$f^{\frak b}_{\frak a}(M)^{n}:=\inf\{f_{{\frak a}R_{\frak p}}^{{\frak b}R_{\frak p}}(M_{\frak p})\,\,|\,\,{\frak p}\in {\rm Spec}(R),\,\,\,\,  \dim R/{\frak p}\geq n\}.$$

\begin{thm}\label{thm2}
Assume that $R$ is a  Noetherian ring  and let  ${\frak b}\subseteq{\frak
a}$ be two ideals of $R$.  Suppose that $M$ is a finitely generated
$R$-module, and let $n$ be a non-negative integer such that the
local cohomology modules $$H_{\frak a}^0(M),\dots, H_{\frak
a}^{t-1}(M)$$ are ${\frak b}$-cofinite, where $t={h_{\frak a}^{\frak
b}(M)}^{n}$. Then ${ f_{\frak a}^{\frak b}(M)}^n=t$. Moreover,
$${ f_{\frak a}^{\frak b}(M)}^1={ h_{\frak a}^{\frak b}(M)}^{1}=\mu_{\frak a}^{\frak b}(M).$$

\proof Since the local cohomology modules $H_{\frak a}^0(M) ,\dots, H_{\frak a}^{t-1}(M)$ are ${\frak b}$-cofinite and ${\frak b}\subseteq{\frak a}$, it follows from  \cite[Lemma 4.2]{HK}
that $H_{\frak a}^0(M) ,\dots, H_{\frak a}^{t-1}(M)$
are ${\frak a}$-cofinite. Hence in view of \cite[Theorem 2.2]{AT} the $R$-module  $\Hom_{R}(R/{\frak a}, H_{\frak a}^{t}(M))$ is finitely generated, and so the set $\Ass_RH^t_{\frak a}(M)$
is finite.  Therefore in view of \cite[Proposition 2.2 and Theorem 2.10]{AN2},   we have
${ f_{\frak a}^{\frak b}(M)}^n={ h_{\frak a}^{\frak b}(M)}^{n}$.

 Now, we show that ${ h_{\frak a}^{\frak b}(M)}^{1}=\mu_{\frak a}^{\frak b}(M).$ To do this, let $i$ be an arbitrary non-negative integer such that $i<\mu_{\frak a}^{\frak b}(M)$.
 Then there is a non-negative integer $t$ such that ${\frak b}^{t}H_{\frak a}^{i}(M)$ is minimax. Hence ${\frak b}^{t}H_{\frak a}^{i}(M)$ is in dimension $<1$, and
 so $h_{\frak a}^{\frak b}(M)^1\geq\mu_{\frak a}^{\frak b}(M)$. To get the equality, let $j$ be an arbitrary non-negative integer such that $j< h_{\frak a}^{\frak b}(M)^1$.
  Then there exists a non-negative integer $t$ such that ${\frak b}^{t}H_{\frak a}^{j}(M)$ is in dimension $<1$.  Therefore, in view of  \cite[Proposition 3.9]{AN2},
  there exists a non-negative  integer $s$ such that $\dim \Supp({\frak b}^{s}H_{\frak a}^{j}(M))<1$, and so $\Supp({\frak b}^{s}H_{\frak a}^{j}(M))\subseteq\Max(R)$. Hence,
  as $\Hom_{R}(R/{\frak a},{\frak b}^{s}H_{\frak a}^{j}(M))$ is finitely generated, it follows that $\Hom_{R}(R/{\frak a},{\frak b}^{s}H_{\frak a}^{j}(M))$ is Artinian.  Since
${\frak b}^{s}H_{\frak a}^{j}(M)$ is ${\frak a}$-torsion, it yields
from  Melkersson's result \cite [Theorem 1.3]{Me1} that ${\frak
b}^{s}H_{\frak a}^j(M)$ is  Artinian. Hence ${\frak b}^{s}H_{\frak
a}^j(M)$ is minimax, and so $\mu_{\frak a}^{\frak b}(M)\geq h_{\frak
a}^{\frak b}(M)^1$. This completes
 the proof of theorem.  \qed \\
\end{thm}

We recall that  an $R$-module $L$ is called {\it skinny or weakly Laskerian}, if each of its homomorphic images has only finitely many associated primes (cf. \cite{Ro} and \cite{DM}).
\begin{cor}\label{cor5*}
Suppose that $R$ is a  Noetherian ring  and let  ${\frak b}\subseteq{\frak a}$ be two ideals of $R$.  Assume that $M$ is a finitely generated $R$-module and let $n$ be a
 non-negative integer such that the local cohomology modules $$H_{\frak a}^0(M) ,\dots, H_{\frak a}^{t-1}(M)$$ are   weakly Laskerian,  where $t={h_{\frak a}^{\frak
b}(M)}^{n}$. Then
${ f_{\frak a}^{\frak b}(M)}^n=t$. Moreover,
$${ f_{\frak a}^{\frak b}(M)}^1={ h_{\frak a}^{\frak b}(M)}^{1}=\mu_{\frak a}^{\frak b}(M).$$
\proof The assertion follows from \cite [Proposition 2.7]{BNS2} and the proof of Theorem \ref{thm2}. \qed \\
\end{cor}
\begin{cor}\label{cor5**}
    Let $R$ be a  Noetherian ring and let  ${\frak b}\subseteq{\frak a}$ be two ideals of $R$. Let $M$ be a finitely generated $R$-module and let $n$ be a non-negative integer
    such that  $\dim\Supp(H_{\frak a}^i(M))\leqslant1$ for all $i<{ h_{\frak a}^{\frak b}(M)}^{n}$. Then
    ${ f_{\frak a}^{\frak b}(M)}^n={ h_{\frak a}^{\frak b}(M)}^{n}$. Moreover,
    $${ f_{\frak a}^{\frak b}(M)}^1={ h_{\frak a}^{\frak b}(M)}^{1}=\mu_{\frak a}^{\frak b}(M).$$
\proof The assertion follows from \cite [Corollary 3.4]{BNS1} and  the proof of Theorem \ref{thm2}. \qed \\
\end{cor}
\begin{cor}\label{cor5***}
Let $R$ be a  Noetherian ring and let  ${\frak b}\subseteq{\frak a}$ be two ideals of $R$. Let $M$ be a finitely generated $R$-module and let $n$ be a non-negative integer
 such that  $\Supp (H_{\frak a}^i(M))\subseteq \Max (R)$ for all $i<{ h_{\frak a}^{\frak b}(M)}^{n}$. Then
${ f_{\frak a}^{\frak b}(M)}^n={ h_{\frak a}^{\frak b}(M)}^{n}$. Moreover,
$${ f_{\frak a}^{\frak b}(M)}^1={ h_{\frak a}^{\frak b}(M)}^{1}=\mu_{\frak a}^{\frak b}(M).$$
\proof The assertion follows from  Corollary 2.7. \qed \\
\end{cor}
\begin{cor}\label{cor6}
Let $R$ be a  Noetherian ring, $M$ a finitely generated $R$-module and  ${\frak a}$  an ideal of $R$ with  $\dim M/\frak aM \leq1$.  Let ${\frak b}$ be a second ideal of $R$
such that ${\frak b}\subseteq{\frak a}$. Then, for every non-negative integer $n$,  we have
 ${ f_{\frak a}^{\frak b}(M)}^n={ h_{\frak a}^{\frak b}(M)}^{n}$. Moreover,
$${ f_{\frak a}^{\frak b}(M)}^1={ h_{\frak a}^{\frak b}(M)}^{1}=\mu_{\frak a}^{\frak b}(M).$$
\proof The assertion follows from Corollary 2.7 and the  fact that, for all $i\geq0$,  $\Supp(H^i_{\frak a}(M))\subseteq \Supp(M/\frak aM)$. \qed \\
\end{cor}

\begin{cor}\label{cor7}
Let $R$ be a  Noetherian ring, ${\frak a}$  an ideal of $R$, and  $M$ a finitely generated $R$-module such that  $\dim M \leq2$.  Let ${\frak b}$ be a second ideal of $R$
such that ${\frak b}\subseteq{\frak a}$. Then, for any non-negative integer $n$,  we have
${ f_{\frak a}^{\frak b}(M)}^n={ h_{\frak a}^{\frak b}(M)}^{n}$. Moreover,
$${ f_{\frak a}^{\frak b}(M)}^1={ h_{\frak a}^{\frak b}(M)}^{1}=\mu_{\frak a}^{\frak b}(M).$$

\proof The assertion follows from  \cite[Corollary 5.2]{CGH} and the proof of Theorem \ref{thm2}. \qed \\
\end{cor}
\begin{defn}
Let $R$ be a commutative Noetherian ring and let $r$ be a positive integer. For any non-negative integer $n$, we say that the  {\it Faltings' local-global principle for the
 in dimension $<$ n of local cohomology modules} holds at level $r$ (over the
ring $R$) if,  for every choice of ideals ${\frak a}$, ${\frak b}$ of
    $R$  and for every choice of finitely generated $R$-module  $M$, it
    is the case that $$ h_{{\frak a}R_{\p}}^{{\frak b}R_{\p}}(M_{\p})^{n}>r
    \text{ for all } {\p}\in \Spec(R) \Longleftrightarrow h_{{\frak a}}^{{\frak b}}(M)^{n}>r.$$
\end{defn}

The following theorem plays a key role in the proof of the main result of this section.
\begin{thm}\label{thm3}
Suppose that $R$ is a  Noetherian ring  and let  $\frak a, {\frak b}$ be two ideals of $R$ such that ${\frak b}\subseteq{\frak a}$.  Assume that $M$ is a finitely generated $R$-module and
 let  $r$ be a positive integer such that  the local cohomology
modules $H_{\frak a}^0(M),\dots, H_{\frak a}^{r-1}(M)$ are ${\frak b}$-cofinite. Then for any non-negative integer $n$,
$$ h_{{\frak a}R_{\p}}^{{\frak b}R_{\p}}(M_{\p})^{n}>r \text{ for all } {\p}\in \Spec(R) \Longleftrightarrow h_{{\frak a}}^{{\frak b}}(M)^{n}>r.$$
\end{thm}

\proof Let $ h_{{\frak a}R_{\p}}^{{\frak b}R_{\p}}(M_{\p})^{n}>r \text{ for all } {\p}\in \Spec(R)$ and let  $i$ be an arbitrary  non-negative integer such that $i\leqslant r$. It
is sufficient for us to show that there is a non-negative integer $t_0$ such
that ${\frak b}^{t_0}H_{\frak a}^i(M)$ is in dimension $<n$. To do this, as $H_{\frak a}^0(M),\dots, H_{\frak a}^{r-1}(M)$ are ${\frak b}$-cofinite and ${\frak b}\subseteq{\frak a}$,
it follows from \cite[Lemma
4.2]{HK} that $H_{\frak a}^0(M),\dots, H_{\frak a}^{r-1}(M)$ are also ${\frak a}$-cofinite, and so in view of  \cite[Theorem
2.2]{AT}, the $R$-module $\Hom_R(R/\frak a, H_{\frak a}^i(M))$  is finitely generated for all $i=0,1,\dots,r$. Hence the set $\Ass_R({\frak b}^tH_{\frak a}^i(M))$ is
finite,  for all  $t\in\mathbb{N}_0$. Thus for all $t\in\mathbb{N}_0$, the set $\Supp ({\frak
    b}^tH_{\frak a}^i(M))$ is a closed subset of $\Spec (R)$ (in the
Zariski topology), and so the descending chain
$$\dots \supseteq \Supp( {\frak b}^tH_{\frak a}^i(M))\supseteq \Supp ({\frak b}^{t+1}H_{\frak a}^i(M))\supseteq \dots$$
is eventually stationary. Therefore there is a non-negative integer $t_0$ such that for each $t\geqslant t_0$, $$\Supp ({\frak b}^tH_{\frak a}^i(M))=\Supp ({\frak b}^{t_0}H_{\frak a}^i(M)).$$
Now, let  ${\p}\in \Spec (R)$ with $\dim R/{\p}\geq n$, and let ${\m}$ be a maximal ideal of $R$ such that ${\p}\subseteq{\m}$. Since $h_{{\frak
a}R_{\m}}^{{\frak b}R_{\m}}(M_{\m})^{n}>r$,  it follows that  there exists an integer
$u\geqslant t_0$ such that $({\frak b}R_{\m})^uH_{{\frak a}R_{\m}}^i(M_{\m})$ is in dimension $<n$. Hence there is a finitely generated submodule $N$ of
 $({\frak b}R_{\m})^uH_{{\frak a}R_{\m}}^i(M_{\m})$ such that $\dim\Supp ({\frak b}R_{\m})^uH_{{\frak a}R_{\m}}^i(M_{\m})/N< n$. Since $\dim R_{\m}/{\frak p}R_{\m}\geq n $, it follows that
$$(({\frak b}R_{\m})^uH_{{\frak a}R_{\m}}^i(M_{\m})/N)_{{\frak p}R_{\m}}=0.$$
Hence the $R_{\m}$-module $(({\frak b}R_{\m})^uH_{{\frak a}R_{\m}}^i(M_{\m}))_{{\frak p}R_{\m}}$ is finitely generated, and so it follows from
 $$(({\frak b}R_{\m})^uH_{{\frak a}R_{\m}}^i(M_{\m}))_{{\frak p}R_{\m}}\cong ({\frak b}^uH_{\frak a}^i(M))_{\p}$$
that $ ({\frak b}^uH_{\frak a}^i(M))_{\p}$ is a finitely generated $R_{\p}$-module for all ${\p}\in \Spec (R)$ with $\dim R/{\frak p}\geq n $. Now,
as $({\frak b}^uH_{\frak a}^i(M))_{\p}$ is ${\frak a}R_{\p}$-torsion, there is
an integer $v\geq1$ such that $({\frak b}^{u+v}H_{\frak a}^i(M))_{\p}=0$, and
so  ${\p}\not\in \Supp ({\frak b}^{t_0}H_{\frak a}^i(M))$. Therefore
$$\Supp ({\frak b}^{t_0}H_{\frak a}^i(M))\subseteq \{ {\p}\in \Spec (R)\,\,|\,\,\dim R/{\frak p}< n\}.$$
Consequently, $\dim\Supp({\frak b}^{t_0}H_{\frak a}^i(M))<n$, for all $i\leq r$, and hence $ h_{{\frak a}}^{{\frak b}}(M)^{n}>r$, as required. \qed \\

The first consequence of Theorem 2.12 is a generalization of the
main result of Raghavan \cite{Ra}.
\begin{cor}\label{cor11}
    The local-global principle (for the in dimension $<n$ of local cohomology modules) holds at level $1$ (over any commutative Noetherian ring).
\end{cor}

\proof The assertion follows from  Theorem \ref{thm3}.\qed
\begin{cor}\label{cor15}
    Let $R$ be a  Noetherian ring  and let  $\frak a, {\frak b}$ be two ideals of $R$ such that ${\frak b}\subseteq{\frak a}$.  Let $M$ be a finitely generated $R$-module and let
     $r$ be a positive integer such that  the local cohomology
    modules $H_{\frak a}^0(M),\dots, H_{\frak a}^{r-1}(M)$ are ${\frak b}$-cofinite. Then
        $$f_{{\frak a}R_{\p}}^{{\frak b}R_{\p}}(M_{\p})>r \,\,\,\,  \text{ for all } {\p}\in \Spec (R) \Longleftrightarrow f_{{\frak a}}^{{\frak b}}(M)>r.$$
\end{cor}

\proof The assertion follows Theorem \ref{thm3} and the fact that $h_{\frak a}^{\frak b}(M)^{0}=f_{\frak a}^{\frak b}(M)$.\qed\\
\begin{cor}\label{cor16*}
    The local-global principle (for the annihilation of local cohomology modules) holds at all levels $r\in \mathbb{N}$, over any (commutative Noetherian) ring $R$ with $\dim R\leq 2$.
\end{cor}

\proof The result follows easily from  \cite[Corollary 5.2]{CGH}, Theorem \ref{thm3},  and the fact that $h_{\frak a}^{\frak b}(M)^{0}=f_{\frak a}^{\frak b}(M)$, for any $R$-module $M$.\qed \\

\begin{cor}\label{cor8}
   Let $R$ be a  Noetherian ring and let  $\frak a,  \frak b$ be two ideals of $R$ such that ${\frak b}\subseteq{\frak a}$.  Let $M$ be a finitely generated $R$-module with $\dim M\leq 2$
    and  $r$  a non-negative integer. Then
    $$h_{{\frak a}R_{\p}}^{{\frak b}R_{\p}}(M_{\p})^{n}>r \,\,\,\,  \text{ for all } {\p}\in \Spec (R) \Longleftrightarrow h_{{\frak a}}^{{\frak b}}(M)^{n}>r.$$
\end{cor}

\proof The result follows easily from  \cite[Corollary 5.2]{CGH} and Theorem \ref{thm3}.\qed \\
\begin{cor}\label{cor9}
    Let $R$ be a  Noetherian ring, $M$ a finitely generated $R$-module and $\frak b \subseteq \frak a$   ideals of $R$. Let $r$ be a positive integer such that
     $\dim\Supp (H_{\frak a}^i(M))\leqslant1$ for all $i<r$. Then
    $$h_{{\frak a}R_{\p}}^{{\frak b}R_{\p}}(M_{\p})^{n}>r \text{ for all } {\p}\in \Spec (R) \Longleftrightarrow h_{{\frak a}}^{{\frak b}}(M)^{n}>r.$$
    In particular,
    $$f_{\frak a}^{\frak b}(M)=\inf\{f_{{\frak a}R_{\frak p}}^{{\frak b}R_{\frak p}}(M_{\frak p})\,\,|\,\,{\frak p}\in {\rm Spec}(R)\}.$$
\end{cor}

\proof The assertion follows from \cite[Corollary 3.4]{BNS1} and the proof of Theorem \ref{thm3}.\qed \\
\begin{cor}\label{cor10}
    Let $R$ be a  Noetherian ring, $M$ a finitely generated $R$-module and  ${\frak a}$  an ideal of $R$ with  $\dim M/\frak aM \leq1$.  Let ${\frak b}$ be a
    second ideal of $R$ such that ${\frak b}\subseteq{\frak a}$. Then, for any positive integer $r$,
    $$h_{{\frak a}R_{\p}}^{{\frak b}R_{\p}}(M_{\p})^{n}>r \text{ for all } {\p}\in \Spec (R) \Longleftrightarrow h_{{\frak a}}^{{\frak b}}(M)^{n}>r.$$
        In particular,
    $$f_{\frak a}^{\frak b}(M)=\inf\{f_{{\frak a}R_{\frak p}}^{{\frak b}R_{\frak p}}(M_{\frak p})\,\,|\,\,{\frak p}\in {\rm Spec}(R)\}.$$
\end{cor}

\proof The assertion follows from \cite[Corollary 3.5]{BNS1} and Theorem \ref{thm3}.\qed \\

Our next corollary is  a generalization  of \cite[Corollary 2.3]{BRS}.
\begin{cor}\label{cor12}
    Let $R$ be a  Noetherian ring and let  $\frak a,  \frak b$ be two ideals of $R$ such that ${\frak b}\subseteq{\frak a}$.  Let $M$ be a finitely generated $R$-module such
    that $\frak a M\neq M$. Then
    $$h_{{\frak a}R_{\p}}^{{\frak b}R_{\p}}(M_{\p})^{n}>\grade_M \frak a\,\,\,\,  \text{ for all } {\p}\in \Spec (R) \Longleftrightarrow h_{{\frak a}}^{{\frak b}}(M)^{n}>\grade_M \frak a.$$
    In particular,
       $$f_{{\frak a}R_{\p}}^{{\frak b}R_{\p}}(M_{\p})>\grade_M \frak a\,\,\,\,  \text{ for all } {\p}\in \Spec (R) \Longleftrightarrow f_{{\frak a}}^{{\frak b}}(M)>\grade_M \frak a.$$
\end{cor}

\proof The assertion follows from the definition of $\grade_M \frak a$ and Theorem \ref{thm3}.\qed \\

The next result is a generalization of Corollary  \ref{cor12}
\begin{cor}\label{cor13}
    Let $R$ be a  Noetherian ring and let  $\frak a,  \frak b$ be two ideals of $R$ such that ${\frak b}\subseteq{\frak a}$.  Let $M$ be a finitely generated $R$-module, and that
     $r\in\{f_{\frak a}(M), f_{\frak a}^1(M), f_{\frak a}^2(M)\}$. Then
    $$h_{{\frak a}R_{\p}}^{{\frak b}R_{\p}}(M_{\p})^{n}>r \,\,\,\,  \text{ for all } {\p}\in \Spec (R) \Longleftrightarrow h_{{\frak a}}^{{\frak b}}(M)^{n}>r.$$
\end{cor}

\proof The assertion follows from  \cite[Theorems  2.3 and 3.2]{BNS1} and  Theorem \ref{thm3}.\qed \\

Before we state the next result recall that for a finitely generated $R$-module  $M$ and for ideals  $\frak a, {\frak b}$ of $R$ with ${\frak b}\subseteq{\frak a}$, the
 ${\frak b}$-cofiniteness dimension $c_{{\frak a}}^{{\frak b}}(M)$ of $M$ relative to $\frak a$ (see \cite{BNS}) is defined by
$${{\rm c}_{\frak a}^{\frak b}(M)}:=\inf\{i\in\mathbb{N}_0\,\,|\,\,
H_{\frak a}^i(M)\text{ is not ${\frak b}$-cofinite}\}.$$
\begin{cor}\label{cor14}
    Let $R$ be a  Noetherian ring and let  $\frak a,  \frak b$ be two ideals of $R$ such that ${\frak b}\subseteq{\frak a}$.  Let $M$ be a finitely generated $R$-module, and that
     $r=c_{{\frak a}}^{{\frak b}}(M) $. Then
    $$h_{{\frak a}R_{\p}}^{{\frak b}R_{\p}}(M_{\p})^{n}>r \,\,\,\,  \text{ for all } {\p}\in \Spec (R) \Longleftrightarrow h_{{\frak a}}^{{\frak b}}(M)^{n}>r.$$
\end{cor}

\proof The assertion follows from the definition of  $c_{{\frak a}}^{{\frak b}}(M)$  and  Theorem \ref{thm3}.\qed \\

The next corollary gives us a short and easy proof of the main result of \cite[Theorem 2.6]{DN1}.

\begin{cor}\label{cor16}
   Suppose that $R$ is a  Noetherian ring and let  $\frak a, {\frak b}$ be two ideals of $R$ such that ${\frak b}\subseteq{\frak a}$.
   Let $M$ be a finitely generated $R$-module and let  $r$ be a positive integer such that  the local cohomology
    modules $H_{\frak a}^0(M),\dots, H_{\frak a}^{r-1}(M)$ are ${\frak a}$-cofinite. Then
    $$\mu_{{\frak a}R_{\p}}^{{\frak b}R_{\p}}(M_{\p})>r \text{ for all } {\p}\in \Spec (R) \Longleftrightarrow \mu_{{\frak a}}^{{\frak b}}(M)>r.$$
\end{cor}

\proof The result follows from the proofs of Theorems \ref{thm2} and \ref{thm3}.\qed \\

We are now ready to state and prove the  main theorem of this section, which shows that Faltings' local-global principle  for the in dimension $<n$
 of local cohomology modules is valid at level $2$ over any commutative Noetherian $R$. This generalizes the main result of Brodmann et al. in \cite[Theorem 2.6]{BRS}.
\begin{thm}\label{thm4}
    The local-global principle (for the  in dimension $<n$  of local cohomology modules) holds over any (commutative Noetherian) ring $R$  at level $2$.
\end{thm}

\proof Let $M$ be a finitely generated $R$-module and let  $\frak a, {\frak b}$ be two ideals of $R$ such that ${\frak b}\subseteq{\frak a}$. We must show that
    $$h_{{\frak a}R_{\p}}^{{\frak b}R_{\p}}(M_{\p})^{n}>2 \,\,\,\,  \text{ for all } {\p}\in \Spec (R) \Longleftrightarrow h_{{\frak a}}^{{\frak b}}(M)^{n}>2.$$
 To do this, it is enough for us to show that, if $h_{{\frak a}R_{\p}}^{{\frak b}R_{\p}}(M_{\p})^n>2 \,\,\,\,  \text{ for all } {\p}\in \Spec R$,
then $h_{{\frak a}}^{{\frak b}}(M)^n>2$. In view of Corollary \ref{cor11}, we need to show that there exists a non-negative integer $u$ such that the $R$-module
${\frak b}^uH_{\frak a}^2(M)$ is in dimension $<n$. Since $h_{{\frak a}R_{\p}}^{{\frak b}R_{\p}}(M_{\p})^{n}>2$, analogous to the proof of Theorem \ref{thm3},
for each ${\p}\in \Spec (R)$ with $\dim R/{\p}\geq n$, there is $t_{\p}\in \Bbb{N}_0$ such that $({\frak b}^{t_{\p}}H_{\frak a}^i(M))_{\p}=0, i=1,2.$ Furthermore,
 there exists a non-negative integer $s$ such that ${\frak b}^sH_{\frak a}^i(\Gamma_{\frak b}(M))=0$ for all $i\geq0$. Now, let  $\bar{M}=M/\Gamma_{\frak b}(M)$.
 Then from the short exact sequence
$$0\To \Gamma_{\frak b}(M) \To M \To \bar{M} \To 0,$$
we obtain the long exact sequence
$$\hspace{12mm}(H_{\frak a}^1(M))_{\p}\To (H_{\frak a}^1(\bar{M}))_{\p}\To (H_{\frak a}^2(\Gamma_{\frak b}(M)))_{\p}\To (H_{\frak a}^2(M))_{\p}\To (H_{\frak a}^2(\bar{M}))_{\p}.\hspace{13mm}(\dag)$$
Hence, it follows from  \cite[Lemma 9.1.1]{BS} that $({\frak b}R_{\p})^{k_{\p}}H_{{\frak a}R_{\p}}^1(\bar{M}_{\p})=0$, for some integer $k_{\p}\in\nat_0$.
 Moreover, by  \cite[Lemma 2.1.1]{BS}, there exists $x\in{\frak b}$ which
 is a non-zerodivisor on $\bar{M}$. Then $x^{k_{\p}}H_{{\frak a}R_{\p}}^1(\bar{M}_{\p})=0$. Now, the short exact sequence
 $$0\To  \bar{M}_{\p}\stackrel{x^{k_{\p}}}\To \bar{M}_{\p}\To \bar{M}_{\p}/ x^{k_{\p}}
 \bar{M}_{\p}\To 0,$$ induces the exact sequence
 $$H_{{\frak a}R_{\p}}^0(\bar{M}_{\p}/ x^{k_{\p}}
 \bar{M}_{\p})\To H_{{\frak a}R_{\p}}^1(\bar{M}_{\p})\stackrel{x^{k_{\p}}}\To H_{{\frak a}R_{\p}}^1(\bar{M}_{\p}).$$
Hence the $R_{\p}$-module $H_{{\frak a}R_{\p}}^1(\bar{M}_{\p})$ is a
  homomorphic image of $H_{{\frak a}R_{\p}}^0(\bar{M}_{\p}/x^{k_{\p}} \bar{M}_{\p})$, and so it is a finitely generated $R_{\p}$-module, for all
  ${\p}\in\Spec (R)$ with $\dim R/{\p}\geq n$.  It therefore follows from
  \cite [Theorem 2.2]{Ho} that $H_{\frak a}^1(\bar{M})$ is in dimension $<n$. Therefore in view of  \cite[Theorem 2.2]{AT}, the $R$-module  $\Hom_R(R/{\frak a}, H_{\frak a}^2(\bar{M}))$
  is also in dimension $<n$, and so by \cite[Lemma 2.6]{MNS} the set $(\Ass_{R}\Hom_R(R/{\frak a}, H_{\frak a}^2(\bar{M})))_{\geq n}$
is finite and consequently the set  $(\Ass_{R} H_{\frak a}^2(\bar{M})))_{\geq n}$ is finite, and so for every non-negative integer $t$, the set
 $(\Ass_{R}{\frak b}^t H_{\frak a}^2(\bar{M})))_{\geq n}$ is also finite. Thus for all $t\in\nat_0$, the set $(\Supp {\frak b}^{t}H_{\frak a}^2(\bar{M}))_{\geq n}$
is a closed subset of $\Spec (R)$ (in the
Zariski topology), and so the descending chain
$$\dots \supseteq (\Supp{\frak b}^tH_{\frak a}^2(\bar{M}))_{\geq n}\supseteq (\Supp {\frak b}^{t+1}H_{\frak a}^2(\bar{M}))_{\geq n}\supseteq \dots$$
is eventually stationary. Therefore there is a non-negative integer $t_0$ such that for each
 $t\geqslant t_0$, $$(\Supp {\frak b}^tH_{\frak a}^2(\bar{M}))_{\geq n}=(\Supp {\frak b}^{t_0}H_{\frak a}^2(\bar{M}))_{\geq n}.$$
Now, as  $({\frak b}^{t_{\p}}H_{\frak a}^2(M))_{\p}=0$, it follows from the exact sequence $(\dag)$ and \cite[Lemma 9.1.1]{BS}
that there is a non-negative integer $v_{\p}\geq t_0$ such that $$(\frak b^{v_{\p}}H_{\frak a}^2(\bar{M}))_{\p}=0.$$
Hence $(\frak b^{t_{0}}H_{\frak a}^2(\bar{M}))_{\p}=0$ for all ${\p}\in \Spec (R)$ with $\dim R/{\p}\geq n$, and so $$\Supp {\frak b}^{t_0}H_{\frak a}^2(\bar{M})
\subseteq \{ {\p}\in \Spec (R)\,|\,\,\dim R/{\frak p}< n\}.$$
Now, let $u:=s+t_0$. Then, it easily follows from the exact sequence  $(\dag)$ and  \cite[Lemma 9.1.1]{BS} that
$$\Supp ({\frak b}^uH_{\frak a}^2(M))\subseteq \Supp({\frak b}^{t_{0}}H_{\frak a}^2(\bar{M})).$$
Consequently, $\Supp {\frak b}^{u}H_{\frak a}^2({M})
\subseteq \{ {\p}\in \Spec (R)\,\,|\,\,\dim R/{\frak p}< n\}$, and so ${\frak b}^{u}H_{\frak a}^2({M})$ is in dimension $<n$, as required.  \qed \\
\begin{cor}\label{cor21}
    The local-global principle (for the in dimension $<n$ of local cohomology modules) holds over any (commutative Noetherian) ring $R$  with $\dim R\leq 3$.
\end{cor}
\proof  The assertion follows from Corollary \ref{cor11}, Theorem \ref{thm4} and \cite[Exercise 7.1.7]{BS}. \qed \\

As a consequence of Theorem \ref{thm4}, the following corollary shows that the local-global principle for the
annihilation of local cohomology modules holds at level $2$ over $R$ and at all levels whenever $\dim R\leq 3$. These  reprove  the main results of Brodmann et al. in \cite{BRS}.
\begin{cor}\label{cor22}
    The local-global principle (for the annihilation of local cohomology modules) holds over any (commutative Noetherian) ring $R$  at level $2$.
\end{cor}
\proof The assertion follows from Theorem \ref{thm4} and the fact that $h_{\frak a}^{\frak b}(M)^{0}=f_{\frak a}^{\frak b}(M)$, where ${\frak b}\subseteq{\frak a}$ are ideals of $R$
and $M$ a finitely generated $R$-module.\qed \\

 \begin{cor}\label{cor23}
    The local-global principle (for the annihilation of local cohomology modules) holds over any (commutative Noetherian) ring $R$  with $\dim R\leq 3$.
 \end{cor}
 \proof  The assertion follows from Corollary \ref{cor21} and the fact that $h_{\frak a}^{\frak b}(M)^{0}=f_{\frak a}^{\frak b}(M)$, where ${\frak b}\subseteq{\frak a}$
 are ideals of $R$ and $M$ a finitely generated $R$-module.\qed \\
\section{\textbf{Local-global principle and associated primes of local cohomology modules}}
It will be shown in this section that the Faltings' local-global principle, for the
in dimension $<n$ of local cohomology modules,  holds at all levels $r\in \mathbb{N}$  whenever the ring $R$ is a homomorphic image of a Noetherian Gorenstein ring. This
generalizes the main result of Raghavan \cite{Ra}.

Also, as a generalization of the main results of Brodmann-Lashgari and Quy, we  prove a finiteness result about associated primes of local cohomology modules.
The main results are Theorems
\ref{thm7} and  \ref{thm6}.

\begin{thm}\label{thm7}
Suppose that $R$ is a  Noetherian ring  which is  a homomorphic image of a Gorenstein ring. Then the local-global principle (for the  in dimension $<n$  of local cohomology modules)
holds at all levels $r\in\mathbb{N}$.
\end{thm}
\proof Let $M$ be a finitely generated $R$-module and let $\frak a,  \frak b$ be two ideals of $R$ such that ${\frak b}\subseteq{\frak a}$, and
that $ h_{{\frak a}R_{\p}}^{{\frak b}R_{\p}}(M_{\p})^{n}>r$ for all ${\p}\in \Spec (R)$. We must show that $h_{{\frak a}}^{{\frak b}}(M)^{n}> r$.
To this end, let ${\p}\in \Spec (R)$ with $\dim R/{\p}\geq n$. There exists a maximal ideal $\m$ of $R$ such that
$\dim R/{\p}=\dim R_{\m}/{\frak p}R_{\m}\geq n$. Moreover, in view of hypothesis for all non-negative integer $i$ with $i\leq r$
 there exists a non-negative integer $t$ such that the $R_{\m}$-module $({\frak b}R_{\m})^tH_{{\frak a}R_{\m}}^i(M_{\m})$ is in dimension $<n$.
 Hence there is a finitely generated submodule $N$ of $({\frak b}R_{\m})^tH_{{\frak a}R_{\m}}^i(M_{\m})$ such that $\dim\Supp ({\frak b}R_{\m})^tH_{{\frak a}R_{\m}}^i(M_{\m})/N< n$,
  and so it follows from  $\dim R_{\m}/{\frak p}R_{\m}\geq n $ that
$$(({\frak b}R_{\m})^tH_{{\frak a}R_{\m}}^i(M_{\m})/N)_{{\frak p}R_{\m}}=0.$$
Consequently, the $R_{\p}$-module  $$(({\frak b}R_{\m})^tH_{{\frak a}R_{\m}}^i(M_{\m}))_{{\frak p}R_{\m}}\cong ({\frak b}^tH_{\frak a}^i(M))_{\p}$$
 is finitely generated for all ${\p}\in \Spec (R)$ with $\dim R/{\frak p}\geq n $. Since the $R_{\p}$-module  $({\frak b}^tH_{\frak a}^i(M))_{\p}$ is ${\frak a}R_{\p}$-torsion,
 it follows that there exists an integer $v\geq1$ such that $({\frak a}R_{\p})^{v}({\frak b}^{t}H_{\frak a}^i(M))_{\p}=0$, and thus $({\frak b}^{t+v}H_{\frak a}^i(M))_{\p}=0$.
  Therefore for all $i\leq r$ and ${\p}\in \Spec (R)$ with $\dim R/{\p}\geq n$, there exists an integer $t_{\p}$ such that  $({\frak b}^{t_{\p}}H_{\frak a}^i(M))_{\p}=0$.
  Hence, for every  ${\p}\in \Spec (R)$ with $\dim R/{\p}\geq n$, we have $ f_{{\frak a}R_{\p}}^{{\frak b}R_{\p}}(M_{\p})>r$. Thus, in view of \cite[Theorem 9.3.7]{BS},
  we have $ \lambda_{{\frak a}R_{\p}}^{{\frak b}R_{\p}}(M_{\p})>r $ for all ${\p}\in \Spec (R)$ with $\dim R/{\p}\geq n$. Consequently $\lambda_{{\frak a}}^{{\frak b}}(M)_{n}>r $,
  and so in view of \cite[Theorem 2.14]{DN2}, $f_{{\frak a}}^{{\frak b}}(M)_{n}>r $. Now, it follows from \cite[Theorem 3.10]{AN2} that $h_{{\frak a}}^{{\frak b}}(M)^{n}>r $, as required.
\qed\\
\begin{thm}\label{thm8}
Assume that $R$ is a  Noetherian ring and let $M$ be a finitely generated $R$-module. Suppose that   $\frak a,  \frak b$ are two ideals of $R$ such that ${\frak b}\subseteq{\frak a}$,
 and let $n$ be a non-negative integer. Then
$$h_{{\frak a}}^{{\frak b}}(M)^{n}\geq\depth({\frak b},M)   \Longleftrightarrow h_{{\frak a}}^{n}(M)\geq\depth({\frak b},M).$$
\end{thm}
\proof Since $h_{{\frak a}}^{n}(M)\leq h_{{\frak a}}^{{\frak b}}(M)^{n}$, it is enough for us to show that if $h_{{\frak a}}^{{\frak b}}(M)^{n}\geq\depth({\frak b},M)$,
then $h_{{\frak a}}^{n}(M)\geq\depth({\frak b},M)$. To do this, let $\depth({\frak b},M)=s$ and we use induction on $s$. For $s=0$ there is nothing to show. So assume that $s>0$
and the result has been proved for $s-1$. Since $\depth({\frak b},M)>0$, it follows that $\frak b$ contains an element $x$ which is a non-zerodiviser on $M$. Moreover,
as $h_{{\frak a}}^{{\frak b}}(M)^{n}\geq s$, there is an integer $t_{0}$ such that, for all $i<s$, the $R$-module ${\frak b}^{t_{0}}H_{\frak a}^i(M)$ is in dimension $<n$.
Hence in view of \cite[Proposition 3.9]{AN2} there exists an integer $t$ such that, for all $i<s$,
$\dim\Supp ({\frak b}^{t}H_{\frak a}^i(M))<n$. In addition, the exact sequence
 $$0\To  {M} \stackrel{x^{t}}\To {M}\To {M}/ x^{t}
 {M}\To 0,$$
induces the long exact sequence
$$\hspace{12mm} H_{\frak a}^{s-2}(M) \stackrel{x^{t}}\To H_{\frak a}^{s-2}({M})\To H_{\frak a}^{s-2}({M}/ x^{t}{M})\To H_{\frak a}^{s-1}(M) \stackrel{x^{t}}\To H_{\frak a}^{s-1}({M}).\hspace{13mm}$$
Now, it follows from \cite[Lemma 2.9]{DN2} that there exists an integer $l$ such that for all $i<s-1$,
$$\dim\Supp ({\frak b}^{l}H_{\frak a}^i({M}/ x^{t}{M}))<n.$$
Therefore, using \cite[Proposition 3.9]{AN2} we see that there exists an integer $v$ such that the $R$-module ${\frak b}^{v}H_{\frak a}^i({M}/ x^{t}
{M})$ is in dimension $<n$ for all $i<s-1$, and thus $h_{{\frak a}}^{{\frak b}}({M}/ x^{t}{M})^{n}\geq s-1$. Consequently in view of the inductive hypothesis
$h_{{\frak a}}^{n}({M}/ x^{t}{M})\geq s-1$, and so the $R$-module $H_{\frak a}^{s-2}({M}/ x^{t}{M})$ is in dimension $<n$. Therefore, it follows from
\cite[Proposition 2.12]{MNS} and the exact sequence
$$0\To  H_{\frak a}^{s-2}({M})/x^{t}H_{\frak a}^{s-2}({M}) \To H_{\frak a}^{s-2}({M}/ x^{t}{M})\To (0:_{H_{\frak a}^{s-1}({M})}x^{t})\To 0,$$
that the $R$-module $ (0:_{H_{\frak a}^{s-1}({M})}x^{t})$   is in dimension $<n$. Now, let ${\p}\in \Spec (R)$ with $\dim R/{\p}\geq n$. Then the $R$-module
$(0:_{(H_{\frak a}^{s-1}({M}))_{\p}}x^{t})$   is finitely generated. On the other hand, as $\dim\Supp ({\frak b}^{t}H_{\frak a}^{s-1}({M}))<n$,
 we obtain that $( {\frak b}^{t}H_{\frak a}^{s-1}({M})_{\p}=0$, and so $ x^{t}(H_{\frak a}^{s-1}({M}))_{\p}=0$. Therefore the $R_{\p}$-module
$$(H_{\frak a}^{s-1}({M}))_{\p}=(0:_{(H_{\frak a}^{s-1}({M}))_{\p}}x^{t})$$ is  finitely generated for all ${\p}\in \Spec (R)$ with $\dim R/{\p}\geq n$.
 Consequently, $f_{{\frak a}}^{n}(M)\geq s$, and so in view of \cite[Theorem 2.10]{MNS},  $h_{{\frak a}}^{n}(M)\geq s$, as required. \qed\\
\begin{thm}\label{thm5}
Let $R$ be a Noetherian ring, $M$  a finitely generated $R$-module, $\frak a$ an ideal of $R$ and  $r$  a positive integer such that the $R$-modules
${\frak a}^tH_{\frak a}^0(M)$,\dots,${\frak a}^tH_{\frak a}^{r-1}(M)$ are in dimension $<2$ for some $t\in\mathbb{N}_0$. Then the $R$-module $\Hom_R(R/{\frak a},H_{\frak a}^r(M))$
is finitely generated and the $R$-modules $H_{\frak a}^0(M)$,\dots,
$H_{\frak a}^{r-1}(M)$ are ${\frak a}$-cofinite. In particular the set $\Ass_R H^r_{\frak a}(M)$ is finite.
\end{thm}

\proof Since the $R$-modules
${\frak a}^tH_{\frak a}^0(M)$,\dots,${\frak a}^tH_{\frak a}^{r-1}(M)$ are in dimension $<2$, it follows that there exist finitely generated submodules $L_i$ of ${\frak a}^tH_{\frak a}^i(M)$
such that    for all $i=0,1,\dots,r-1$,  $\dim \Supp({\frak a}^tH_{\frak a}^i(M)/L_i)<2$. Therefore, for each ${\p}\in \Spec (R)$ with $\dim R/{\p}\geq2$, we have
$$({\frak a}R_{\p})^tH_{{\frak a}R_{\p}}^i(M_{\p})\cong({\frak a}^tH_{\frak a}^i(M))_{{\p}}\cong (L_i)_{\p}.$$
Hence the $R_{\p}$-module $({\frak a}R_{\p})^tH_{{\frak a}R_{\p}}^i(M_{\p})$ is
finitely generated, for all $i=0,1,\dots,r-1$. Now, as $({\frak a}R_{\p})^tH_{{\frak a}R_{\p}}^i(M_{\frak p})$ is ${\frak a}R_{\p}$-torsion,
so there exists a non-negative integer $s$ such that $({\frak a}R_{\p})^{t+s}H_{{\frak a}R_{\p}}^i(M_{\p})=0$.
Therefore, in view of \cite[Proposition 9.1.2]{BS}, the $R_{\p}$-module $H_{{\frak a}R_{\p}}^i(M_{\p})$  is
finitely generated  for every ${\p}\in \Spec (R)$ with $\dim R/{\p}\geq 2$ and for all $i=0,1,\dots,r-1$. It therefore follows from
\cite[Proposition 3.1]{BNS1} that  $\Hom_R(R/{\frak a}, H_{\frak a}^r(M))$ is finitely generated and the
$R$-modules $H_{\frak a}^0(M)$,\dots, $H_{\frak a}^{r-1}(M)$ are ${\frak a}$-cofinite.\qed \\

\begin{cor}
Let $R$ be a Noetherian ring,  $M$ a finitely generated $R$-module,
$\frak a$ an ideal of $R$ and  $r$  a positive integer such that the
$R$-modules $H_{\frak a}^0(M)$,\dots, $H_{\frak a}^{r-1}(M)$ are in
dimension $<2$. Then the $R$-module $\Hom_R(R/{\frak a},H_{\frak
a}^r(M))$ is finitely generated and  the $R$-modules $H_{\frak
a}^0(M)$,\dots, $H_{\frak a}^{r-1}(M)$ are ${\frak a}$-cofinite.
\end{cor}

\proof The assertion follows from Theorem 3.3 and  \cite[Proposition 2.12]{MNS}. \qed \\
\begin{cor}\label{cor24}
    Let $R$ be a Noetherian ring, $M$  a finitely generated $R$-module, $\frak a$ an ideal of $R$ and  $r$  a positive integer such that the $R$-modules
    ${\frak a}^tH_{\frak a}^0(M)$,\dots, ${\frak a}^tH_{\frak a}^{r-1}(M)$ have finite support, for some $t\in\mathbb{N}_0$. Then
     for any  ideal $\frak b$ of $R$ with $\frak b\subseteq \frak a$,
    $$h_{{\frak a}R_{\p}}^{{\frak b}R_{\p}}(M_{\p})^{n}>r\,\,\, \text{ for all } {\p}\in \Spec (R) \Longleftrightarrow h_{{\frak a}}^{{\frak b}}(M)^{n}>r.$$
\end{cor}

\proof Since the set $\Supp ({\frak a}^tH_{\frak a}^i(M))$ is finite, for all $i<r$, it follows that the $R$-module ${\frak a}^tH_{\frak a}^i(M)$
is in dimension $<2$ for all $i<r$. Now  the  assertion  follows from  Theorems 3.3 and 2.12. \qed \\

The following theorem, which is  the second our main result in this section, generalizes the main results of Brodmann-Lashgari \cite{BL} and Quy \cite{Qu}.
\begin{thm}\label{thm6}
    Let $R$ be a Noetherian ring, $M$  a  finitely generated $R$-module,
    and $\frak a$ an ideal of $R$. Let $r$ be a non-negative integer such that $\frak a^tH^{i}_{I}(M)$ is
    in dimension $<2$ for all $i<r$ and for some $t\in\mathbb{N}_0$.  Then, for any minimax submodule $N$ of
    $H^{r}_{\frak a}(M)$,   the $R$-module $\Hom_{R}(R/\frak a, H^{r}_{\frak a}(M)/N)$ is finitely generated. In particular,
    the set $\Ass_{R}(H^{r}_{\frak a}(M)/N)$ is finite.
\end{thm}
\proof In view of Theorem \ref{thm5}, the $R$-module $\Hom_{R}(R/\frak a, H^{r}_{\frak a}(M))$ is
finitely generated and so the $R$-module $0:_{N}\frak a$ is also finitely generated . Hence according to  Melkersson's result
\cite[Proposition 4.3]{Me2}, $N$ is $\frak a$-cofinite, (note that $N$ is minimax). Moreover, the short exact
sequence $$0 \longrightarrow N \longrightarrow H^{r}_{\frak a}(M)
\longrightarrow H^{r}_{\frak a}(M)/N \longrightarrow 0$$ induces the
following exact sequence,$$\Hom_{R}(R/\frak a,H^{r}_{\frak a}(M))
\longrightarrow \Hom_{R}(R/\frak a,H^{r}_{\frak a}(M)/N) \longrightarrow
\Ext^{1}_{R}(R/\frak a,N).$$ Now, as the $R$-modules $\Hom_{R}(R/\frak a,H^{r}_{\frak a}(M))$ and $\Ext^{1}_{R}(R/\frak a,N)$ are finitely generated,
 it follows that the $R$-module $\Hom_{R}(R/\frak a, H^{r}_{\frak a}(M)/N)$ is also finitely generated, as
required.\qed\\

\begin{center}
{\bf Acknowledgments}
\end{center} The authors are deeply grateful to the referee for his/her careful
reading of the paper and valuable suggestions. Also, we would like to thank Professor Kamran Divaani-Aazar for his reading of the first draft and useful  discussions.
 Finally, we would like to thank from the Institute for Research in Fundamental Sciences (IPM), for the  financial support.

\end{document}